%% file: arxiv.tex
\documentclass[11pt]{article}

\usepackage[english]{babel}

\usepackage[a4paper,top=1.5in,bottom=1.4in,left=1.6in,right=1.6in,marginparwidth=1.75cm]{geometry}

\usepackage{amsmath}
\usepackage{graphicx}
\usepackage{amssymb}
\usepackage{bbm}
\usepackage{todonotes}
\usepackage{float}
\usepackage{authblk}
\usepackage{imakeidx}
\usepackage{cite}
\usepackage[colorlinks=true, allcolors=black]{hyperref}
\allowdisplaybreaks
\usepackage[nottoc]{tocbibind}

\let \le \leqslant
\let \ge \geqslant

\usepackage{amsmath}
\usepackage{amssymb}

\numberwithin{equation}{section}

\newtheorem{thm}{Theorem}[section]
\newtheorem{defi}[thm]{Definition}
\newtheorem{prop}[thm]{Proposition}

\newtheorem{cor}[thm]{Corollary}

\input{macro}

\usepackage{cite}
\usepackage{enumitem}

\title{Behind the equivalent descriptions of the \\
Loewner energy}

\author{Yilin Wang\footnote{Department of Mathematics, Eidgen\"ossische Teschnishe Hochschule, 
Z\"urich, Switzerland. Email: \protect{yilin.wang@math.ethz.ch}}}

\begin{document}

\maketitle

\begin{abstract}
This expository paper retells the story behind the paper \emph{Equivalent descriptions of the Loewner energy}, which uncovered the link between random fractal SLE curves and the K\"ahler geometry of universal Teichm\"uller space via the Loewner energy. It is written for the Proceedings of ICBS as the summary paper for the 2026 Frontiers of Science award.
\end{abstract}

\input{core.tex}

\vspace{10pt}

\noindent \textbf{Acknowledgment:} I thank Steffen Rohde, Fredrik Viklund, and Wendelin Werner for their influence, support from the beginning of my career, and comments on an earlier version of this memoir, and Yuliang Shen for helpful discussions while I worked on the award paper.  I am also grateful to the FSA committee and ICBS for their appreciation and for giving me the fun of writing this memoir, and to Viola Giovannini, Dorian Martino, and Catherine Wolfram for their feedback. 

My research is supported by the Swiss State Secretariat for Education, Research and Innovation (SERI) MB25.00004 and the Simons Foundation International [SFI-MPS-PP-00012621-19].

\medskip 

\noindent \textbf{AI declaration:} No AI was used or consulted. 

\bibliographystyle{abbrv}
\bibliography{ref}




\end{document}

%% file: macro.tex
\usepackage[dvipsnames]{xcolor}

\newcommand{\abs}[1]{\left\lvert #1 \right \rvert}

\newcommand{\brac}[1]{\left \langle #1 \right \rangle}

\newcommand{\mc}[1]{\mathcal{#1}}
\newcommand{\m}[1]{\mathbb{#1}}

\def\ie{i.e.,\,}

\renewcommand\Im{\operatorname{Im}}

\def\SLE{\operatorname{SLE}}

\def\mob{\mathrm{M\ddot{o}b}}
\def\Diff{\operatorname{Diff}}
\def\QS{\operatorname{QS}}

\def\g{\gamma}

\def\D{\Delta}

\def\t{\theta}

\def\s{\sigma}

\def\o{\omega}
\def\O{\Omega}
\def\vare{\varepsilon}

\def\Chat{\widehat{\m{C}}}

 \def\Id{\operatorname{Id}}

\def\dd{\mathrm{d}}
\def\vol{\mathrm{vol}}
\def\ii{\mathfrak{i}}
\newcommand{\ee}{\mathrm{e}}

\newcommand{\ad}[1]{\overline{#1}}

\def\detz{\mathrm{det}_{\zeta}}

\def\1{\mathbf{1}}


%% file: core.tex
\section{Introduction}
I had the great pleasure of receiving the Frontier of Sciences Award for my paper \emph{Equivalent descriptions of the Loewner energy} \cite{W2} at the 2026 International Congress of Basic Science (ICBS) in Beijing. I wasn't expecting it, and I was even more pleased by how much I enjoyed the conference, especially the talks by other awardees. Certainly, the topics were broad, featured award-worthy breakthroughs, and were delivered by curious, interesting minds. But other than that, the unusual format of ICBS also contributed to its uniqueness: speakers only talked about one paper instead of a plethora of impressive results; the winning paper could be quite old, so the talks were filled with mature hindsight; and it was one of the works that we were most proud of and to which we had a great emotional attachment. 

So instead of writing a review article, I decided to write a different kind of article to tell the story behind this paper, which aligns better with why I enjoyed ICBS. I wrote an expository article for Notices of the AMS \cite{wang2023nams} on the same paper a couple of years ago, which is another reason I want to write something different rather than merely repeating it.

\subsection{The reversibility of chordal Loewner energy} 

The story started with a question Steffen Rohde emailed me after seeing my first paper \cite{W1} showing the reversibility of chordal Loewner energy. Let me explain the background. 

Consider the upper half-plane $\m H :=  \{z \in \m C \mid \Im (z) > 0\}$ and a simple curve $\gamma$ in $\m H$ going from $0$ to $\infty$, we call such a curve a \emph{chord}. The Loewner transform (introduced by Charles Loewner in 1923 \cite{Loewner23}) encodes $\gamma$ into a driving function $W : \m R_+ \to \m R$ in the following way. We first parametrize $\gamma$ by capacity. This means that $\gamma(0) = 0$ and for $t \ge 0$, the unique conformal map (biholomorphic function) $g_t : \m H \smallsetminus \gamma (0,t] \to \m H$ that has expansion at $\infty$ as $g_t (z) = z + O(1/z)$, has the next term being exactly $2t/z$, namely
$$g_t (z)  = z + \frac{2t}{z} + o \left(\frac1z\right),$$
and the coefficient $2t$ is the \emph{half-plane capacity} of $\gamma (0,t]$ which can be thought as a measurement of its size.

Intuitively, the family of conformal maps $(g_t)_{t \ge 0}$ \emph{unzip} progressively the curve $\gamma$ from the bottom, and sends the tip of the slit, $\gamma (t)$, to a point $W (t):= g_t (\gamma (t)) \in \m R$. The function $W : \m R_+ \to \m R$ is a continuous function, and is called the \emph{Loewner driving function} of the chord $\gamma$. This encoding is very effective: we can recover the family $(g_t)_{t \ge 0}$, and hence the curve $\gamma$, from $W$ as it satisfies a simple ordinary differential equation 
$$
\begin{cases} \partial_t g_t (z) = \dfrac{2}{g_t (z) - W_t}, \qquad \forall z \in \m H\\
g_0 (z) = z  
\end{cases}
$$
as long as it is defined.

Oded Schramm \cite{Schramm2000} used this method to describe a family of random chords by requiring their driving function to be a multiple of Brownian motion, giving rise to Schramm--Loewner Evolution SLE, which describes the fractal scaling limits of interfaces in critical statistical mechanics models, e.g., \cite{LSW04LERWUST,SS_harmonic,CDHKS}. SLE curves now play a central role in two-dimensional random conformal geometry, statistical mechanics, and a probabilistic approach to quantum gravity and conformal field theory. More precisely, if $W = \sqrt \kappa B$ where $\kappa \ge 0$, then the random curve in $\m H$ is called the \emph{chordal SLE$_\kappa$} (there are some subtleties in showing that they are actually curves \cite{Rohde_Schramm}).  

Let's come back to deterministic chords. The \emph{chordal Loewner energy} of a deterministic chord $\gamma$ is defined \cite{W1,FS15,Dub_comm} as the Dirichlet energy of its driving function:
$$I_{\m H; 0, \infty} (\gamma) := \frac12 \int_0^\infty (W'(t))^2 \, \dd t$$
if $W$ is absolutely continuous (in particular $W'$ is well-defined), and $\infty$ otherwise. 
It is easy to check that if $\gamma$ is scaled by a factor $\lambda > 0$, then its driving function is changed to $t \mapsto \lambda W(\lambda^{-2} t)$ which has the same Dirichlet energy. In particular, it allows us to define, for any chord in a simply connected domain $D$ going from a prime end (boundary point) $a$, to another different prime end $b$ as
$$I_{D; a, b} (\gamma) := I_{\m H; 0, \infty} (\psi (\gamma))$$
where $\psi$ is a conformal map $D \to \m H$ sending respectively $a,b$ to $0,\infty$.

This definition was inspired by the study of large deviations of $\SLE$, namely, the rarity of the event where $\SLE_\kappa$ stays close to $\gamma$ as $\kappa \to 0$. Exploring the relation to $\SLE$, and in particular, the known reversibility of $\SLE$ \cite{Zhan}, I gave a probabilistic proof of the reversibility of the Loewner energy.
\begin{thm}[W.~\cite{W1}]
    For any chord $\gamma$ in $(D;a,b)$, we have 
    $$I_{D; a, b} (\gamma) = I_{D; b, a} (\gamma).$$
    In other words, if $\gamma$ is a chord in $\m H$ going from $0$ to $\infty$, then
    $$I_{\m H; 0, \infty} (\gamma) = I_{\m H; 0, \infty} (-1/\gamma),$$
    where $-1/\gamma$ is the image of $\gamma$ under the map $\iota: \m H \to \m H: z\mapsto -1/z$.
\end{thm}
A priori, the statement has nothing to do with random curves. But it would be surprising for someone who doesn't know SLE, since the driving function $W(t)$ of $\gamma$ is determined by the part $\gamma(0,t]$; reversing the curve (by the map $\iota$) breaks the dependence. In fact, we don't know how the driving function $\widetilde W$ of $-1/\gamma$ is explicitly related to that of $\gamma$, and it is likely the wrong question to ask. However, this theorem shows that $W$ and $\widetilde W$ have the same Dirichlet energy. 

Nowadays, this theorem has several proofs without using SLE, one of them follows from the award paper, but see also \cite{sung_rev}.

\subsection{Birth of the loop Loewner energy}
In the email Steffen wrote me, he seemed greatly intrigued by this reversibility and was convinced that the Loewner energy would be fun to work with for a complex analyst like him. He pointed out that the Loewner energy can be extended to a Jordan curve (a simple closed curve) in $\Chat = \m C \cup \{\infty\}$.

Take a continuously parametrized Jordan curve $\gamma : [0,1] \to \Chat$ with \emph{root} $\gamma (0) = \gamma (1)$. We define the \emph{loop Loewner energy} 
$$I^L (\gamma; \gamma (0)) := \lim_{\vare \to 0+} I_{\Chat \setminus \gamma [0,\vare]; \gamma (\vare), \gamma (0)} (\gamma [\vare, 1]) \in [0,\infty].$$
It is a generalization of the chordal energy, since for a chord $\eta$ in $\Chat \smallsetminus \m R_+$ connecting $0$ and $\infty$ (we no longer care about the orientation), 
\begin{equation}
    \label{eq:loop_chord}
I^L (\eta \cup \m R_+; \infty) = I_{\Chat \smallsetminus \m R_+; 0, \infty} (\eta).
\end{equation}

Steffen asked me, could it be true, that the loop energy does not depend on the choice of the root $\gamma (0)$?

I was immediately absorbed in this question, because if the answer is positive, then the Loewner energy has much more symmetry than the reversibility. On the other hand, it would mean, in the example in \eqref{eq:loop_chord}, $\m R_+$ that used to be the boundary of the domain $\Chat \smallsetminus \m R_+$ for a chordal energy, plays the same role as part of a curve for which we are computing the energy. It sounded too good to be true.
I then tried to check some examples and soon realized that if $\gamma$ is a concatenation of circular arcs, the loop Loewner energy indeed does not depend on the root. It is quite a big family of curves already, so I wrote back to Steffen that I was convinced that the root-invariance is \emph{true}. 

I then went to Seattle for a semester to work on this problem with him. We introduced the loop Loewner energy on paper and proved its root-invariance:
\begin{thm}[Rohde--W. \cite{RW}]
    The loop Loewner energy of a Jordan curve $\gamma$ does not depend on its root. It is invariant under all M\"obius transformations of $\Chat$ and equals $0$ if and only if $\gamma$ is a circle.
\end{thm}
M\"obius transformations are the conformal automorphisms of $\Chat$.  This theorem says that, in a certain sense, the Loewner energy measures the roundness of a Jordan curve in a conformally invariant way. 

The Loewner energy of Jordan curves has later inspired a stream of work connecting to many other fields of mathematics that I will discuss later. But the techniques and new concepts coming from the proof also turn out to be fruitful. 
We had to prove an equivalence between the H\"older regularity of $\gamma$ and that of its driving function, and use an approximation by a family of \emph{piecewise geodesic} curves that we also introduced, which has also inspired many further works \cite{SWW_shears,MRW2026,BJRW2025,Wang_optimization,Bonk_Eremenko,peltola_wang}.

At that time, the root-invariance of the loop Loewner energy was quite surprising, as it showed a lot of symmetry (by changing the root along the curve) that was not apparent from the definition, even more mysterious than the reversibility. It is less surprising now thanks to the \emph{Equivalent descriptions}, or the relations to the SLE loop measures discovered later \cite{BJ2,CarfagniniWang2024}.

 \section{Equivalent descriptions}
 After returning to Zurich, I couldn't help but feel even more puzzled by the root invariance despite the proof we had, and I set out to find an equivalent definition of loop Loewner energy that clearly does not depend on the root. There were no conjectures to prove, nor did I have any evidence that such an equivalent definition existed, only the belief that the world would only make sense if it existed. 
 
 \subsection{Main result}

After many failed attempts, I finally found it:

\begin{thm}[Theorem~1.4 of \cite{W2}]\label{thm:main}
    Let $\gamma$ be a Jordan curve in $\m C$, $\O$ and $\O^*$ be respectively the bounded and unbounded connected components of $\m C \smallsetminus \gamma$. Let $f : \m D \to \O$ and $g : \m D^* \to \O^*$ be two conformal maps such that $g (\infty) = \infty$. Then,
    \begin{equation}
        \label{eq:loewner_main}
    I^L (\gamma) = \frac{1}{\pi}\int_{\m D} \left| \frac{f''}{f'} (z)\right|^2  |\dd z|^2+ \frac{1}{\pi}\int_{\m D^*} \left|\frac{g''}{g'} (z)\right|^2 |\dd z|^2 + 4\, \log\frac{ |f'(0)| }{|g'(\infty)|}.
    \end{equation}
\end{thm}
Here $\m D = \{z \in \m C \mid |z| < 1\}$ and $\m D^* = \{z \in \m C \mid |z| > 1\}$.

In particular, there is no special marked point on the curve involved in the expression, so the root-invariance and reversibility follow automatically. But, more surprisingly and well beyond a better explanation it provides, the right-hand side is exactly the \emph{universal Liouville action} introduced by Takhtajan and Teo \cite{TT06} when they studied the K\"ahler geometry of the universal Teichm\"uller space, to which I will come back later.

If the curve passes through $\infty$, there is a similar formula:
\begin{prop}[Theorem~1.1 of \cite{W2}]\label{prop:infinite}
    If $\infty \in \gamma$, then we let $H$ and $H'$ denote the two connected components of $\mathbb C \smallsetminus \Gamma$, $\mathbb H$ and $\mathbb H^*$ be the upper and lower half-planes, and $f \colon \mathbb H \to H$ and $g\colon \mathbb H^* \to H^*$ be the conformal maps that fix $\infty$. The Loewner energy equals in this case
\begin{equation}
    \label{eq:Loewner_main_2}
I^L(\gamma) = \frac{1}{\pi}\int_{\mathbb H} \left| \frac{f''}{f'} (z)\right|^2  |\dd z|^2+ \frac{1}{\pi}\int_{\mathbb H^*} \left|\frac{g''}{g'} (z)\right|^2 |\dd z|^2.
\end{equation}
\end{prop}

It is hard to describe the euphoria of finally convincing myself of these expressions and the urge to let everyone know about it. It was very memorable.

\subsection{Outline of the proof}
I will try to describe the thought process that finally led to the expressions and their proof. There is no room to discuss the details of the other failed attempts. In my mind, it felt like a long walk through a dark, foggy landscape, until it became crystal clear as the sun rose and pierced the darkness. 
But it was also a walk guarded by my advisor, Wendelin Werner, to whom I am infinitely indebted. 

At first, I came to him with a wild guess, and he told me if he thought it was probably a dead end. Then, a week later, I came up with something else and started over again, and again. Every week I learned something new, dreamed the wildest dream, and reality-checked with him. Wendelin later described it as watching ``a child trying to climb on a rock and never giving up --- who keeps taking new fresh runs at it full of positive energy and enthusiasm,'' and that he was ``looking forward to hear about the new attempts every week.'' 
Some of the missteps may be obvious to avoid by an expert like him, but for the novice I was, it was nevertheless a precious experience to map my own mental landscape. 

\subsubsection{Infinite curves}
As already mentioned, the Loewner energy was introduced as the large deviation rate function of SLE. There were celebrated and insightful results \cite{Dub_couplings,Quantum_zipper,MoT} about how some versions of SLE curves are related to some random surfaces built from Gaussian free fields (GFF) on $\m C$. I remember the shock when I heard Scott Sheffield's talk about how two random trees can be mated to produce a random sphere where the branches of the trees become SLE curves. It has stuck with me. However, at that time, there was not yet a clean link between the \emph{loop} Loewner energy and the SLE loop measure, nor were the SLE/GFF couplings explicitly known in this loop case (now it is known \cite{AHS_quantum_disk}). 

However, after taking a few leaps of faith, I arrived at the statement of Proposition~\ref{prop:infinite}.  
The reason for looking at the curves passing through $\infty$ first is that they relate to the chordal energy more easily, such as the ones in \eqref{eq:loop_chord}. The right-hand side of \eqref{eq:Loewner_main_2} is thought of as the large deviation rate function of GFF.  

Once we have the guess, actually showing it does not need to involve SLE or GFF. The most natural step is to first check whether the right-hand side of \eqref{eq:Loewner_main_2} satisfies the same properties as the Loewner energy.  The crucial property is the additivity, which is nothing more than the observation that 
$$\frac12\int_0^\infty (W'(t))^2 \, \dd t =\frac12 \int_0^s (W'(t))^2 \, \dd t + \frac12\int_s^\infty (W'(t))^2 \, \dd t$$
for all $s \ge 0$ in the case of the Loewner energy. The additivity of the Loewner energy translates to an additive identity between integrals on the right-hand side of \eqref{eq:Loewner_main_2}. This was the most technical part of the proof, showing that this additive identity indeed holds.

Once we have the additivity, we use it to check some examples of curves $\g$ to test the identity. The most accessible ones are those curves in \eqref{eq:loop_chord} where $\eta$ has a piecewise linear driving function. Finally, we conclude by approximating a curve $\gamma$ passing through $\infty$ by a sequence of $\eta_n \cup [n, \infty)$, where $\eta_n$ is a chord in $(\m C \smallsetminus [n, \infty); n,\infty)$ with piecewise linear driving function.

\medskip

Proposition~\ref{prop:infinite} is already a major step, but $\infty$ seems to still be a special point on the curve, and we need to move $\infty$ away from the curve to handle bounded curves. Meanwhile, the right-hand side of \eqref{eq:Loewner_main_2} looked so clean that I thought it would be surprising that nobody had looked at it. Then I found the monograph \cite{TT06} on universal Teichm\"uller space and the expression of the universal Liouville action \eqref{eq:loewner_main}, which looked very similar. 

The monograph \cite{TT06} became \emph{the} math book that I read the most, although back then, it seemed difficult to penetrate. I carried it with me when I traveled back to Shanghai for Christmas in 2017. Luckily, I found out that Yuliang Shen from Soochow University, a one-hour train ride from Shanghai, had worked extensively on related topics. I reached out and was really fortunate to be received by him. I showed up in his office covered with snow, and showed him my formula in Proposition~\ref{prop:infinite}. He confirmed that he believed it should be exactly the universal Liouville action. This gave me some confidence to try to show it. 

\subsubsection{Determinants of Laplacians}

As the Loewner energy is M\"obius invariant, there should be a direct way to show that \eqref{eq:loewner_main} and \eqref{eq:Loewner_main_2} are equal, when we map a curve that passes through $\infty$ to a bounded one using a M\"obius transformation. It was not straightforward with my limited compute by hand though. However, I am glad I didn't brute-force it through. 

The alternative is to find a more conceptual way to move $\infty$ away from the curve. To make $\infty$ like any other point, we consider the curve $\gamma$ in $\widehat {\m C} \simeq \m S^2$ instead of in $\m C$. Inspired by Dub\'edat's work \cite{Dub_couplings} relating SLE to determinants of Laplacians, the Loewner energy must also be related to determinants of Laplacians. 

 Indeed, we equip the sphere $\m S^2$ with a Riemannian metric $g = \ee^{2 \s} g_0$ that is Weyl-equivalent to the round metric $g_0$ (the metric induced from $\m S^2 \subset \m R^3$). 
Let $\g \subset \m S^2$ be a \emph{smooth} Jordan curve dividing $\m S^2$ into two connected components $D_1$ and $D_2$. 
Denote by $\D_{D_i,g}$ the Laplace--Beltrami operator with Dirichlet boundary condition on $(D_i, g)$.
The regularity assumption on $\gamma$ is such that the zeta-regularized determinants of Laplacians $\detz (\D)$ below are well-defined. 
I introduced the \emph{Loewner potential} $\mc H (\cdot, g)$ on the space of smooth Jordan curves:
 $$ \mc H (\g, g) : = \log \detz'(\D_{\m S^2, g}) -  \log \vol_g(\m S^2)  - \log \detz (\D_{D_1, g}) -  \log \detz (\D_{D_2,g}). $$     

Starting from Proposition~\ref{prop:infinite}, we obtain the following result.
\begin{thm}[Theorem 7.3 of {\cite{W2}}] \label{thm_energy_determinant}
We have:
\begin{enumerate}
\item The functional $\mc H$ is conformally invariant, \ie $\mc H (\cdot, g) = \mc H (\cdot, g_0)$;
\item \label{item_energy_determinant}  Let $\g$ be a smooth Jordan curve on $\m S^2$. 
We have the identity
$$I^L(\g) = 12 \mc H(\g, g) - 12 \mc H(\mc C, g)  = 12 \log \frac{\detz(\D_{\m D_1, g})\detz(\D_{\m D_2, g})}{\detz(\D_{D_1, g}) \detz(\D_{D_2, g})},$$
where $\mc C$ is any circle and $\m D_1$ and $\m D_2$ are the two components of  $\widehat {\m C} \smallsetminus \mc C$. 
\end{enumerate}
\end{thm}

The key input in the proof is the Polyakov--Alvarez formula \cite{OPS} which computes the ratio between $\detz(\D_{\m D_i, g})$ and $\detz(\D_{D_i, g})$.

This theorem bridges Proposition~\ref{prop:infinite} and Theorem~\ref{thm:main} since $\infty$ is no longer a special point on the sphere and can be ``rotated away" from $\gamma$.  
So for a smooth Jordan curve that does not pass through $\infty$, we use Theorem~\ref{thm_energy_determinant} and the Polyakov--Alvarez formula again to deduce the identity in Theorem~\ref{thm:main}. For a Jordan curve with less regularity, we use an approximation argument by smooth Jordan curves to complete the proof of Theorem~\ref{thm:main}.

\medskip 

Although Loewner potential was just a tool in this proof, its generalizations have been explored further in other contexts, such as for multichordal SLE, which led to a new proof of the Shapiro conjecture on the enumeration of rational functions and the derivation of classical BPZ equations \cite{peltola_wang,alberts_pole}, and for multiply-connected domains, very recent relations to Koebe uniformization theorem \cite{Luo_Maibach,luo2026}.

\subsection{Weil--Petersson Teichm\"uller space}

\subsubsection{Weil--Petersson curves}
Theorem~\ref{thm:main} led to an immediate answer about which curves have finite Loewner energy. There are three terms in \eqref{eq:loewner_main}. The last one is always finite, and the finiteness of the first two terms is, in fact, equivalent for Jordan curves. So $I^L(\gamma) < \infty$ if and only if $f''/f' \in L^2 (\m D)$. The latter is a condition on the regularity of the curve. 
\begin{defi}[Takhtajan--Teo]
    A Jordan curve in $\m C$ is called a \emph{Weil--Petersson curve} if $f''/f' \in L^2 (\m D)$. 
\end{defi}
Weil--Petersson curves are also quasicircles, so they are sometimes called Weil--Petersson quasicircles. Theorem~\ref{thm:main} then implies:
\begin{cor}
    The Loewner energy of $\gamma$ is finite if and only if $ \gamma$ is a Weil--Petersson curve.
\end{cor}
There are now more than thirty equivalent definitions of Weil--Petersson curves, involving Schwarzian derivatives of $f$, spectral properties of Grunsky operators, planar Coulomb gas, Brownian loop measures, renormalized area of minimal surfaces in the hyperbolic $3$-space $\m H^3$, renormalized volume in $\m H^3$, $\beta$-numbers from the traveling salesman problem, random geometry SLE, M\"obius energy from knot theory, etc. Some of those equivalent definitions are summarized and derived in \cite{B2025}, but see also \cite{johansson_viklund,ViklundWang2024,BBVW,cui00,michelat2021loewner}. 
Gehring and Hag wrote the book \cite{gehring2012ubiquitous} on \emph{The ubiquitous quasidisk}, the same should be said about Weil--Petersson curves. 

\subsubsection{Conformal welding and universal Teichm\"uller space}
Let me describe the context where Takhtajan and Teo introduced the universal Liouville action, the right-hand side of \eqref{eq:loewner_main}.  

We first identify a Jordan curve $\g$ with a homeomorphism of the unit circle $\m S^1 = \partial \m D = \partial \m D^*$ as follows. 
By Carath\'eodory's theorem, the conformal map $f : \m D \to \O$  extends continuously to a homeomorphism between the closures $\ad {\m D} \to \ad {\O}$. Similarly for $g : \m D^* \to \O^*$. In particular, $f$ and $g$ restricted to $\m S^1$ define two homeomorphisms $\m S^1 \to \g$. The \emph{welding homeomorphism} is the circle homeomorphism $\varphi : = g^{-1} \circ f|_{\m S^1}$ which compares these two homeomorphisms. To fix the ambiguity in the choices of $f$ and $g$, we assume that $\g$ passes through $-1,\ii,1$ (which is possible after a M\"obius transformation) and $f, g$ both fix these three points. In other words, we obtain the welding homeomorphism of $\g$ in
$$ \{\varphi \text{ is a homeomorphism } \m S^1 \to \m S^1 \mid \varphi \text{ fixes } -1, \ii, 1\} =\mob (\m S^1)  \backslash \operatorname{Hom} (\m S^1)$$
where $\mob (\m S^1)$ is the group of M\"obius transformations preserving $\m S^1$.

The converse operation consists of finding a Jordan curve $\g$ passing through $-1,\ii,1$ and the corresponding normalized conformal maps $f$ and $g$, whose welding homeomorphism $g^{-1} \circ f|_{\m S^1}$ is a given circle homeomorphism $\varphi$. 
A solution may not exist, and if one exists, it may not be unique. 
The complete characterization of circle homeomorphisms for which there is a unique solution to conformal welding is a hard open question. 

However, a classical result in quasiconformal mappings says that if the circle homeomorphism is \emph{quasisymmetric}, then the solution to the conformal welding problem exists and is unique.
We use $\QS(\m S^1)$ to denote the group of quasisymmetric circle homeomorphisms, which consists of $\varphi \in \operatorname{Hom}(\m S^1)$ such that
$$\exists M \ge 1, \, \forall \t \in \m R, \, \forall t \in (0,\pi), \,  \frac1M \le \abs{\frac{\varphi(e^{\ii (\t + t)}) -\varphi( e^{\ii \t})}{\varphi(e^{\ii \t} )- \varphi(e^{\ii (\t-t)})}} \le M.$$   The corresponding welded Jordan curves are called \emph{quasicircles}. 
The homogeneous space 
$$T(1) : = \mob (\m S^1) \backslash \QS(\m S^1)$$ is the \emph{universal Teichm\"uller space} which has an infinite-dimensional complex Banach manifold structure. For those who are familiar with Riemann surfaces, the Teichm\"uller space (the space of marked complex structures) of any surface of negative Euler characteristic embeds holomorphically into $T(1)$, hence ``universal."

We will only consider quasicircles and quasisymmetric circle homeomorphisms. So it is safe to think that conformal welding provides a correspondence between Jordan curves and circle homeomorphisms. However, conceptually they are very different. For instance, the group structure on $\QS(\m S^1)$ allows us to compose two circle homeomorphisms, and therefore it makes sense to speak about a homogeneous space. But it is not so obvious why one would want to ``compose" two Jordan curves and what does that even mean.

However, this shift in viewpoint allows us to consider geometric structures on the space of quasicircles and describe the tangent bundle, Riemannian metrics, symplectic forms, etc. 
In particular, physicists already wondered if it is possible to put a \emph{homogeneous K\"ahler metric} on the subspace $\mob (\m S^1) \backslash \Diff^\infty (\m S^1)$ since it is the space of reparametrizations of a closed string (and they only care about smooth reparametrizations)  \cite{BowickRajeev1987string,witten88}. A \emph{homogeneous K\"ahler metric} is a compatible symplectic form $\o (\cdot, \cdot)$ and a Riemannian metric $\brac{\cdot, \cdot}$ that are invariant under the right-composition and compatible with the complex structure (encoded in an operator $J$ on the tangent bundle such that $J^2 = -\Id$, which plays the role of ``multiplication by $\ii$''). Being compatible means $\o (\cdot, J (\cdot)) = \brac{\cdot, \cdot}$.

\medskip

Bowick and Rajeev \cite{BowickRajeev1987string} showed that there is indeed a homogeneous K\"ahler metric, and only one (up to a scaling factor)!

\medskip

The completion of $\mob (S^1)\backslash \Diff^\infty(S^1)$ for this unique metric is denoted as 
$T_0(1)$  
and called \emph{Weil--Petersson Teichm\"uller space} \cite{TT06}. 
It is further identified explicitly as 
$$T_0(1) = \mob (\m S^1)  \backslash \{\varphi \mid \varphi \text{ is absolutely continuous and } \log |\varphi'| \in H^{1/2} (\m S^1) \}$$
by Yuliang Shen \cite{shen13}.
Moreover, $\varphi \in T_0(1)$ if and only if the welded curve is a Weil--Petersson quasicircle \cite{TT06}, or equivalently, a curve with finite Loewner energy by Theorem~\ref{thm:main}. 

The link between Loewner energy and this K\"ahler metric goes further.
One special feature of K\"ahler metrics is that they always admit a locally defined K\"ahler potential. Here, the Weil--Petersson metric has a globally defined K\"ahler potential:
 
\begin{thm}[Takhtajan--Teo {\cite[Cor.\,II.4.2]{TT06}}]
    The universal Liouville action ${\mathbf S} 
    : T_0(1) \to \m R_+$ is a K\"ahler potential for the Weil--Petersson metric. That is, up to a positive scaling factor,
    $$\ii \partial \bar\partial {\mathbf S} = \o, $$
    where $\partial$ and $\bar \partial$ are the $(1,0)$ and $(0,1)$ components of de Rham differential on the complex manifold $T_0(1)$.
\end{thm}
This shows that this unique K\"ahler metric on $T_0(1)$ is completely encoded into the universal Liouville action, hence, the Loewner energy by Theorem~\ref{thm:main}.

\subsubsection{A step back}
At this point, I was totally in awe. How can the Loewner energy that arises from the distribution of random curves from statistical mechanics have anything to do with the unique K\"ahler metric on this infinite-dimensional space of curves? I find myself parsing more and more on the \emph{Why} than the \emph{How}; there is not yet a completely satisfactory answer. The journey started from trying to understand the \emph{Why} behind the root-invariance, but it only led to a bigger mystery. Perhaps the journey of wonder never ends.

\subsubsection{Beyond the bridge}
Nevertheless, a bridge is now built between two fields of mathematics, random conformal geometry and Teichm\"uller theory, that were previously remote. It opens a goldmine of new mathematics. Let me list a few works that have benefited from this bridge other than those already mentioned: 

\begin{itemize}[left = 2 pt, topsep= 5pt, itemsep=2pt]
\item \cite{W3}: Since SLE are tightly connected to Brownian loop measure \cite{LawlerWerner2004} and Werner's loop measure \cite{Werner_loop}, this work expressed the Loewner energy and characterizes  Weil--Petersson curves in terms of those loop measures.
    \item \cite{VW1}: Inspired by SLE/GFF couplings, including imaginary geometry \cite{IG1} and the quantum zipper \cite{Quantum_zipper}, this work revealed the relations between the Loewner energy and Dirichlet energy of functions defined on complex domains mimicking the couplings.
    \item \cite{ViklundWang2024}: Inspired by the SLE duality, we obtained results about foliations of $\Chat$ by Weil--Petersson curves and the duality between the foliation energy (Loewner--Kufarev energy) and Loewner energy, which in turn suggests the reversibility of the whole-plane SLE for all $\kappa > 8$ (the case $\kappa \le 8$ was known to be true by \cite{zhan_rev_whole,IG4}).
    \item \cite{ang232}: Inspired by the reversibility of foliation energy, this work proved the conjecture about the reversibility of the whole-plane SLE for all $\kappa > 8$. 
    \item \cite{GQW2025,BJ1}: The geometric viewpoint on the Loewner energy consists of deforming the curve infinitesimally at all points, which is very different from the Loewner evolution viewpoint, which explores along the curve. This global geometric viewpoint underlies the construction of Virasoro representation from SLE loop measures, and the characterization of SLE loop measures as the unique Malliavin--Kontsevich--Suhov measures.
      \item \cite{sung2023quasiconformal}: This work uses the Loewner evolution to give a new proof of the variational formula of the Loewner energy/universal Liouville action. 
    \item \cite{TBBT}: This work further explores the variational formula of the Loewner energy in terms of the Brownian loop measure.
  \item \cite{fan2025QI,BJ2}: These works explore the group structures in SLE loop measures, and derive their Cameron--Martin theorem with respect to composition of weldings.
  \item \cite{FVW2026}: This work aims to provide an expression for the Loewner energy directly from the welding, without recourse to solving the welding problem. 
  \item \cite{krusell2024rho}: This work generalized the results of the award paper to the large deviation rate function of $\SLE_\kappa(\rho)$ curves.
    \item \cite{BBVW,martino2026WP}: Renormalized volume of convex co-compact hyperbolic $3$-manifolds provide a K\"ahler potential of Weil--Petersson metrics on Teichm\"uller spaces of surfaces.  Inspired by it, we obtained a new expression of the Loewner energy as a renormalized volume of $\m H^3$ as the ``universal analog.'' This should open a new venue of studying the holographic principle of random curves SLE.
\end{itemize}
Since I dislike packing a plethora of results into a talk, I will refrain from saying more here. They are just for those who would love to read more.

\section{Outlook}

If you ask me, what part of this journey do I cherish the most? I would say the time before the dawn of finding the \emph{Equivalent descriptions}. It was a period of my PhD during which I did not produce much. But I learned to face infinite possibilities and the reality of not understanding, and to sit with the darkness; I had to draw my own map, seek help from experts and literature, and fight through it all with uncertainty. Although Wendelin was supportive and hoping to see me solve the problem, he also suggested that I look into other things in parallel so that I would not be discouraged by being stuck. I just couldn't. The question haunted me day and night. No matter what other directions I looked, it came back to me.  The moment it came to light, everything made sense. 

\medskip 

Eight years have passed, I still believe, the journey of wonder never ends.